\documentclass[12pt]{article}

\usepackage{amsmath, amsfonts, amssymb, amsthm} \parskip1pt

\def\cO{ {\cal O}}

\def\bN{{\mathbb N}}

\def\bR{{\mathbb  R}}

\def\bx{{\bf x}}

\def\by{{\bf y}}

\def\bz{{\bf z}}
\def\bZ{{\mathbb Z}}

  \def\pro{\noindent {\bf{Proof. }}}

 \def\Card{{\rm Card}}

\def\build#1_#2^#3{\mathrel{\mathop{\kern 0pt#1}\limits_{#2}^{#3}}}

\def\smallsquare{\vbox{\hrule\hbox{\vrule height 1 ex\kern 1
ex\vrule}\hrule}}
\def\cqfd{\hfill \smallsquare\vskip 3mm}

\def\cB{{\cal B}}

\def\cE{{\cal E}}

\def\vol{{\rm vol}}

\newtheorem{thm}{Theorem}
\newtheorem{lem}{Lemma}
\newtheorem*{prop}{Proposition}
\newtheorem*{prob}{Problem}

\newtheorem*{Definition}{Definition}

\date{ }
\begin{document}

\title{ Inhomogeneous approximation with coprime integers and lattice orbits}
\author{Michel Laurent \;\; $\&$ \;\; Arnaldo Nogueira}
\maketitle


\noindent  ABSTRACT -- 
Let  $(\xi,y)$  be a  point in $\bR^2$ and  $\psi : \bN \rightarrow \bR^+$   a  function.  
We investigate the problem of the existence of infinitely many pairs $p, q$ of coprime integers such that
$$
| q \xi +p -y | \le \psi( | q|).
$$
We prove that the existence is ensured  for every point $(\xi,y)$ whose coordinate   $\xi$ is irrational  when   
$
\psi(| q|) = c| q|^{-1/2} 
$ 
for  some suitable coefficient
$
c= c(\xi,y).
$ 
 An optimal metrical statement  is also established. We link the subject with density exponents
 of lattice orbits in $\bR^2$.

\vskip 9mm

\section{ Introduction and results}
Minkowski has proved that for every  real irrational number $\xi$ and every real number $y$ not belonging to 
$\bZ \xi + \bZ$, there exist infinitely many pairs of integers $p, q$ such that
$$
 | q \xi + p - y | \le { 1 \over 4  | q| }. \leqno{(1)}
$$
See for instance Theorem II in Chapter 3 of Cassels' monograph \cite{Cas}. The statement is optimal in the sense that the approximating function
$\ell \mapsto  (4 \ell  )^{-1}$ cannot be decreased. Note that the restriction $y \notin \bZ \xi +\bZ$ can  be dropped at the cost of replacing 
the upper bound $ (4| q|)^{-1}$
by  $ c | q|^{-1}$ for any constant $c$ greater than $1/\sqrt{5}$. When $y=0$, the primitive point 
$({p\over \gcd(p,q)},{q\over \gcd(p,q)})$ remains a solution to (1), therefore we  may   moreover require  that 
the pair of integers  $p,q$ satisfying (1) be coprime.
 However, for a non-zero real number 
$y$, this extra requirement is far from being obvious to satisfy. We obtain  such a result  with a weaker accuracy 
  of order $1/\sqrt{|  q |}$. 

\begin{thm}
Let $\xi$ be an irrational real number and   $y$  a real number. There exist infinitely many pairs of coprime integers $(p,q)$ 
such that
$$
| q \xi + p - y | \le { c  \over   \sqrt{| q |}} \quad { \rm with} \quad c = 2 \sqrt{3} \max(1, | \xi|)^{1/2}| y|^{1/2}. \leqno{(2)}
$$
\end{thm}

Theorem 1 will be deduced in Section 2 from our results \cite{LN}  of effective density   for $SL(2,\bZ)$-orbits in $\bR^2$.
However, it will clearly appear that our method for proving Theorem 1 provides more than necessary. 
  Theorem 1 may probably be improved. We address the following 
  
\begin{prob}
 Can we replace the approximating function $\psi(\ell) = c \, \ell^{-1/2}$ occurring in Theorem 1
 by a smaller one, possibly $\psi(\ell ) = c\, \ell^{-1}$ ? 
\end{prob}

We shall further discuss this  problem in Section 4 for   the function $\psi(\ell) = 2 \, \ell^{-1}$, offering  some hints  
and indicating the difficulties which then arise. 
It turns out that the  approximating function $\psi(\ell) = \ell^{-1}$ is permitted  for almost all pairs $(\xi,y)$
of real numbers relatively to  Lebesgue measure   and 
that $-1$ is the critical exponent. The last assertion follows from  the following metrical statement:

\begin{thm}
Let $\psi : \bN \mapsto \bR^+$ be a function.  Assume that $\psi$  is non-increasing, tends to $0$ at infinity and   that
for every  positive integer  $c$ there exist positive real numbers $c_1<c_2$ satisfying 
$$
c_1 \psi(\ell) \le \psi(c \, \ell) \le c_2  \psi(\ell) ,  \quad \forall  \ell \ge 1. \leqno{(3)}  
$$
Assume furthermore    that
$$
\sum_{\ell \ge 1} \psi(\ell) = + \infty.
$$
Then,  for almost all pairs $(\xi,y)$ of real numbers  there exist infinitely  many primitive points  $(p,q)$ 
such that
$$
q \ge 1 \quad {\rm and}\quad | q \xi + p - y | \le \psi( q).  \leqno{(4)}
$$
If  $
\sum_{\ell \ge 1} \psi(\ell) 
$
converges, the  pairs $(\xi,y)$  satisfying $(4)$ for infinitely many primitive points  $(p,q)$ form a set of   null
Lebesgue measure.
\end{thm}

 Note that $(4)$ gives  an  additional information on the sign of $q$, and that 
we could  have equivalently  required that $q$ be negative. Such a refinement could 
as well be achieved in the frame of Theorem 1,  with a weaker approximating function of the 
form $\psi(\ell) = \ell^{-\mu}$ for any given  real number $\mu <1/3$, by employing alternatively
Theorem 5 in  Section  9 of \cite{LN}. We  leave the details of proof  to the interested reader, arguing as in Section 2.
 For questions of density involving  signs, see  also    \cite{DaNo}.

 The proof of Theorem 2 is  given in Section 3. It 
combines   standard tools  from  metrical number theory
  with the  ergodic  properties  of the linear action of $SL(2,\bZ)$ on $\bR^2$ \cite{Mo}.
  We refer to  Harman's book \cite{Har} for closely related results. See also the recent overview \cite{BBDV} and 
  the monographs \cite{Sc}, \cite{Sp}. 
  
  Theorem  2 is a metrical statement about  pairs $(\xi,y)$ of real numbers. A natural question is  to understand what happens on  
  each fiber when we fix either $\xi$ or $y$.  In this direction, here is a  partial result
  which will be deduced  from the explicit construction displayed in Section 4.
  
   \begin{thm}
  Let  $\xi$ be an irrational number and let $(p_k/q_k)_{k\ge 0}$ be the sequence of  its convergents. Assume that the series 
  $$
   \sum_{k\ge 0} { 1\over \max(1,\log q_k)} \leqno{(5)}
$$ 
  diverges. Then for almost every real  number  $y$ there exist infinitely many primitive points $(p,q)$ 
  satisfying
  $$
   | q \xi  + p -y | \le { 2 \over | q |  }.
  $$
   Moreover the  series $(5)$ diverges for almost every  real number $\xi$.
  \end{thm}

  We now turn to the second part of the paper devoted to density exponents for lattice orbits in $\bR^2$.
   As already mentioned,
   the approximating function  $\psi(\ell) = c  \,\ell^{-1/2}$ occurring in Theorem 1 is directly connected 
  to the density exponent $1/2$ for $SL(2,\bZ)$-orbits. We intend  to show that this exponent $1/2$ is best possible
  in general. Thus  our method  for proving Theorem 1 
fails   with an approximating function $\ell^{-\mu}$ for $\mu >1/2$. 

We work in the more general setting of  {\it lattices} 
 $\Gamma$ in $SL(2,\bR)$.  Recall that a    lattice $\Gamma$ in $SL(2,\bR)$ is a discrete subgroup for which the  quotient  $ \Gamma \backslash SL(2,\bR)$ has  finite Haar measure.
 We view $\bR^2$ as a space of column vectors on which 
the group of matrices $\Gamma$ acts by left multiplication. 
We  equip $\bR^2$ with   the supremum  norm $|  \, \,\,  |$, and for any matrix $\gamma \in \Gamma$,  we 
 denote as well by $|\gamma |$
 the maximum of the absolute values of the entries of $\gamma$.
 Let us first give a
% some definitions.  
%  A   lattice $\Gamma$ in $SL(2,\bR)$ is a discrete subgroup whose covolume
%for the Haar measure on $SL(2,\bR)$  is finite. Hedlund has shown that any orbit $\Gamma\bx$ in $\bR^2$ is either dense or discrete.
%When $\Gamma = SL(2,\bZ)$, the orbit $\Gamma\bx$ is dense if and only if  the slope of $\bx$ is an irrational number. More information 
% can be found in  Chapter 5 of \cite{Dal}.

\begin{Definition}
Let $\bx$ and $\by$ be two points in $\bR^2$. We denote by $\mu_\Gamma(\bx,\by)$ the supremum, possibly infinite,
 of the exponents $\mu$ such that the inequality
$$
| \gamma\bx - \by | \le | \gamma|^{-\mu} \leqno{(6)}
$$
has infinitely many solutions $\gamma \in \Gamma$.
\end{Definition}

Note that for a fixed $\bx \in \bR^2$, the function $\by \mapsto \mu_\Gamma(\bx,\by)$ is $\Gamma$-invariant. By the ergodicity of the action 
of $\Gamma$ on $\bR^2$, see \cite{Mo}, this function is therefore constant almost everywhere on $\bR^2$. 
We denote by $\mu_\Gamma(\bx) $ its  generic  value and we call $\mu_\Gamma(\bx)$
 the {\it generic density exponent} of the orbit $\Gamma\bx$.

\begin{thm}
The upper bound 
$
\mu_{\Gamma}(\bx) \le {1/ 2}
$
holds true for any point $\bx\in \bR^2$ such that the orbit $\Gamma \bx$ is dense in $\bR^2$. 
\end{thm}

In an equivalent way, Theorem 4 asserts that the upper bound 
$ \mu(\bx,\by) \le 1/2$
holds   for almost all points $\by \in \bR^2$. This bound was  already known  in the case of the unimodular group 
$\Gamma =  SL(2,\bZ)$ as a consequence of  Theorem 3
in \cite{LN}.
 
One may optimistically conjecture that $\mu_{\Gamma}(\bx) = 1/2$ for every   
point $\bx$ such that $\Gamma \bx$ is dense in $\bR^2$, or at least for almost every     point $\bx \in \bR^2$.
In this direction, it follows from  \cite{LN} that the  lower bound 
$$
\mu_{SL(2,\bZ)} (\bx) \ge {1\over 3}
$$
holds  for all points $\bx$ in $\bR^2\setminus\{ {\bf 0} \}$ with irrational slope. 
Weaker lower bounds can as well  be deduced from \cite{MaWe}  which are valid for any lattice $\Gamma \subset SL(2,\bR)$.
Note that the function 
$\bx \mapsto \mu_\Gamma(\bx)$ is
 $\Gamma$-invariant since the quantity  $\mu_\Gamma(\bx)$ obviously depends only on the orbit $\Gamma\bx$. 
 Thus, the generic density exponent $\mu_\Gamma(\bx)$ takes  the same value  for almost all  points $\bx \in \bR^2$.

  \section{ Proof of Theorem 1}
  We first state   a  result  obtained  in  \cite{LN}. In this section, we denote by $\Gamma $ the  lattice   ${\rm SL(2,\bZ)}$.
    For any point  $\bx =  \left(\begin
    {matrix} 
     x_1\cr x_2
      \end{matrix}
      \right)$  
      in $ \bR^2$ with irrational slope $ x_1/x_2$, the 
  orbit $\Gamma \bx$ is  dense in $\bR^2$.   We have obtained in \cite{LN}
effective results concerning  the density of such an orbit. In particular, our estimates are essentially optimal when the target point 
$\by$ has rational slope.

\begin{lem} 
Let  $\bx$ be a  point in $\bR^2$ with irrational slope and  $\by =  \left(\begin{matrix}
 y\cr y\cr
 \end{matrix}
 \right)$
 a point on the diagonal  with $y\not= 0$. 
Then, there exist infinitely many  matrices $\gamma \in \Gamma$ such that
$$
| \gamma \bx - \by | \le {c \over | \gamma |^{1/2}} \quad {\rm with} \quad  c =   2\sqrt{3} | \bx|^{1/2} | y |^{1/2}.\leqno{(7)}
$$
\end{lem}

\pro The point $\by$  has rational slope $1$. Apply Theorem 1 (ii) of \cite{LN} with $a=b=1$.
\cqfd

Put $\bx =  \left(\begin{matrix}
 \xi\cr 1\cr
 \end{matrix}
 \right)$. The point $\bx$ has irrational slope $\xi$ so that Lemma 1 may be applied.
Write  
$
\gamma =  \left(\begin{matrix}
 q_1 &  p_1\cr  q_2& p_2\cr
 \end{matrix}
 \right)
$
 a   matrix provided  by Lemma 1. Then,  the inequality   (7)   gives
 $$
 \begin{aligned}
 \max\left( | q_1 \xi + p_1 -y | ,  | q_2 \xi + p_2 - y | \right) &  \le {c \over \max ( | p_1|, | p_2| , | q_1|, |p_2| )^{1/2}}
 \\
 & \le {c \over \max (  | q_1|, |q_2| )^{1/2}}.
 \end{aligned}
 $$
Therefore, both points $(p,q) = (p_1,q_1)$ and $(p,q) = (p_2,q_2)$ satisfy (2), and since 
the determinant $q_1p_2 - q_2 p_1 = 1$, the two integer points  $(p_1, q_1) $ and $(p_2,q_2)$ are primitive.
 As  there exist  infinitely many matrices $\gamma$ verifying (7), we thus find  infinitely many coprime solutions
to (2). 

\bigskip
\noindent
{\bf Remark.}
The choice of another fixed rational direction would give the same kind of estimate, with  a different constant $c$.
Notice however that our approach  gives more than required, since we get two solutions of (2) forming a matrix $\gamma$ of determinant $1$.
 A better understanding of the  shrinking target problem for  the dense orbit $\Gamma \bx$, not to  a  point $\by$ as in \cite{LN}
but to  a  line in $\bR^2$, may  possibly  lead to  the expected exponent $-1$.

\section{ Proof of Theorem 2}
It is convenient to view the pairs  $(\xi,y)$ occurring in Theorem 2 as    column  vectors  $ \left( \begin{matrix}
\xi \cr y\cr \end{matrix}
\right)$ in $\bR^2$.
 We are  concerned with the set
$\cE(\psi)$ of  vectors   $ \left( \begin{matrix}\xi \cr y\cr\end{matrix}
\right) \in \bR^2$ 
for which  there exist infinitely many primitive integer points $(p,q)$  such that
$$
q \ge 1 \quad {\rm and}\quad  | q\xi +p -y | \le \psi(  q ) . \leqno{(8)}
 $$
For fixed $p,q$,  denote by $\cE_{p,q}(\psi)$ the strip
$$
\cE_{p,q}(\psi) :=\Big \{ \left(\begin{matrix}\xi \cr y\cr\end{matrix}
\right) \in \bR^2 ; \quad  | q\xi +p -y | \le \psi(  q )\Big\}, 
$$
and for every  positive  integer $q$,  let 
$$
\cE_q(\psi) : = \bigcup_{ { p\in \bZ\atop \gcd(p,q) =1}} \cE_{p,q}(\psi)
$$
be the union of all  relevant strips involved  in (8) for fixed $q$. Without loss of generality, we  shall 
assume that $\psi( q) \le 1/2$,  so that the above union is disjoint.
Then $\cE(\psi)$ is equal to  the $ \limsup$ set
$$
\cE( \psi) = \bigcap_{Q \ge 1} \bigcup_{ q  \ge Q} \cE_q(\psi).
$$

As usual when dealing with $\limsup$ set in metrical theory, we  first estimate Lebesgue   measure of pairwise intersections
 of the subsets $\cE_q(\psi),  q  \ge 1$. We establish next    a new kind of zero-one law.

\subsection
{ Measuring intersections}
In this section, we restrict our attention to points located in the unit square $[0,1]^2$.
We denote  by $\varphi$ the Euler  totient function and by $\lambda$ the Lebesgue  measure on $\bR^2$.

\begin{lem}
Let $\psi : \bN \rightarrow [0,1/2]$ be a function. \\
{\rm (i)} For every  positive  integer $q$, we have
$$
\lambda (\cE_q(\psi) \cap [0,1]^2) = { 2 \varphi( q ) \psi( q) \over  q } . 
$$
{\rm (ii)} Let $q$ and $s$ be distinct positive integers. Then,    we have the upper bound
$$
\lambda (\cE_q(\psi) \cap \cE_s(\psi) \cap [0,1]^2)  \le 4 \psi(  q   )\psi(  s ). 
$$
\end{lem}

\pro
Denote by $\chi_q$ the characteristic function of the  interval  
$
\left[-\psi(  q ),  \psi(  q ) \right].
$ 
Then the characteristic function $\chi_{\cE_q(\psi)}$ of the subset $\cE_q(\psi) \subset \bR^2$ is equal to
$$
\chi_{\cE_q(\psi)} (\xi, y) = \sum_{{ p\in \bZ \atop \gcd(p,q) = 1}} \chi_q(q \xi +p - y)= \sum_{{ p\in \bZ \atop \gcd(p,q) = 1}} \chi_q(q \xi -p - y).
$$
Observe that if $\left(\begin{matrix}\xi\cr y\end{matrix}\right)$ belongs to $[0,1]^2$, the indices $p$ of  
 non-vanishing  terms occurring in the last  sum  are located  in the interval $-1 \le p \le q$.  Integrating first with respect to $x$,  we find
$$
\begin{aligned}
\lambda (\cE_q(\psi) \cap [0,1]^2)  & = \int_0^1\int_0^1 \chi_{\cE_q(\psi)} (x, y) dxdy
\\
& =  \sum_{{ p\in \bZ \atop  -1 \le p \le q, \, \gcd(p,q) = 1 }}\int_0^1\int_0^1 \chi_q(q x - p - y)dxdy 
\\
& =
 \int_{1-\psi(q)}^{1 }{-1+y +\psi(q)\over q} dy 
 + \sum_{{1\le  p \le q-2 \atop \gcd(p,q)=1}} \int_0^1{2 \psi(q)\over q}  dy
 \\
 &   \qquad \qquad + \int_0^{1-\psi(q)}{2\psi(q)\over q} dy +  
 \int_{1-\psi(q)}^1{1 -y + \psi(q)\over q}dy
 \\
& = { 2 \varphi( q ) \psi( q) \over  q  }.
\end{aligned}
$$
The  first term appearing  in  the third equality of the above formula 
corresponds to the summation index $p=-1$ and the two last ones  to $p=q-1$.
We have thus  proved (i). 

For the second assertion, we majorize  
$$
\begin{aligned}
\lambda (\cE_q(\psi) \cap \cE_s(\psi)&  \cap [0,1]^2)  = \int_0^1\int_0^1 \chi_{\cE_q(\psi)} (x, y)\chi_{\cE_s(\psi)} (x, y) dxdy
\\
& \le  \int_0^1\int_0^1  \left(\sum_{ p\in \bZ } \chi_q(q x +p - y) \right) \left(\sum_{ r\in \bZ } \chi_s(s x +r - y)\right) dxdy
\\
& = \int_0^1\int_0^1\chi_q(\| q x  - y\|) \chi_s(\| s x  - y \|) dx dy ,
\end{aligned}
$$
where  $\| .  \|$ stands as usual for the distance to the nearest integer. Now,  (ii) follows from the probabilistic independence formula
$$ 
\int_0^1\int_0^1\chi_q(\| q x  - y\|) \chi_s(\| s x  - y \|) dx dy =  4 \psi(  q   )\psi(  s ), 
$$
obtained  by Cassels on page 124 of \cite {Cas} (see  Proof $(ii)$). \cqfd

\subsection
{ A zero-one law}
We say that a subset of $\bR^2$ is a {\it null}  set  if it  has Lebesgue measure 0. 
 A set whose complementary is a null set  is called a {\it full}  set.
The goal of this section is to prove the

\begin{prop}
Let $\psi $ be an approximating function as in Theorem 2. Then the subset $\cE(\psi)$ is either a null set or a full set.
\end{prop}

For proving the proposition, it is convenient to introduce the larger subset
$$
\cE'(\psi ) = \bigcup_{k\ge 1} \cE(k\psi).
$$
In other words, $\cE'(\psi)$ is the set of all points $\left(\begin{matrix}
\xi \cr y\cr
\end{matrix}
\right)$
 in $\bR^2
$ for which there exist  a positive real number 
$\kappa $, depending possibly on $ \left(\begin{matrix}
\xi \cr y\cr
\end{matrix}
\right) $, and infinitely many primitive points $(p,q)$ satisfying
$$
q\ge 1 \quad {\rm and}\quad | q\xi +p -y | \le \kappa \psi(q). \leqno{(9)}
$$
Observe that $\cE(k\psi) \subseteq \cE(k' \psi)$ if $1 \le k\le k'$. In particular, $\cE(\psi)$ is contained in $\cE'(\psi)$. 

\begin{lem}
 Let   $\psi : \bN \rightarrow \bR^+$  be a function tending   to zero at infinity. Then  
the difference  $ \cE'(\psi) \setminus \cE(\psi)$ is  a set of  null Lebesgue measure.
\end{lem}

\pro
We show that all  sets  $ \cE(k\psi),  k\ge 1, $  have the same Lebesgue measure. 
For every real number $y$, denote by $\cE(\psi, y)\subseteq \bR$ the section of $\cE(\psi)$ on the horizontal line  $\bR \times \{ y\}$, i.e.
$$
\cE(\psi, y) = \left\{ \xi \in \bR \, ;  \,   \left(\begin{matrix} \xi \cr  y \end{matrix}\right) \in \cE(\psi)  \right\}.
$$
Then, using $(8)$, we can express
$$
\cE(\psi,y) = \bigcap_{Q\ge 1}\bigcup_{q\ge Q}\bigcup_{{p\in \bZ\atop\gcd(p,q)=1}}\left[ {-p+ y -\psi(q)\over q},{-p+ y +\psi(q)\over q}\right]
$$
as a limsup set of intervals. If we restrict to a bounded part of $\cE(\psi,y)$, the above union over $p$ reduces to a finite one.
Observe  that the centers ${-p+y\over q}$ of these  intervals do not depend on $\psi$, and that their length is multiplied by the constant
factor $k$ when replacing $\psi$ by $k\psi$.  Appealing now to a result due to Cassels \cite{CasB},  we   infer
   that all $\limsup$ sets $\cE(k\psi, y), k\ge 1,$ have the same Lebesgue measure. 
See also  Corollary of Lemma 2.1 on page 30 of
\cite{Har}. 
Notice  that for fixed $k$, the length
${2k\psi(q)\over q}$ of the intervals $\left[ {-p+ y -k\psi(q)\over q},{-p+ y +k\psi(q)\over q}\right]$
 tend to $0$  as  $q$ tends to infinity, as required by  Lemma 2.1.  
 By Fubini, the fibered sets
$$
\cE(k\psi ) = \coprod_{y\in\bR} \Big( \cE(k\psi,y)\times \{ y\}\Big), \quad k\ge 1,
$$
have as well   the same Lebesgue measure in $\bR^2$. \cqfd

\begin{lem}
  Let   $\psi : \bN \rightarrow \bR^+$  be a  function  satisfying the  conditions $(3)$. 
Then  $\cE'(\psi)$ is either a null or a full set.
\end{lem}

\pro 
It is based on the following observation. Let $\left(\begin{matrix}\xi \cr y\cr\end{matrix}\right)$ belong to $\cE'(\psi)$ and let 
$\gamma = \left(\begin{matrix} a& b\cr c& d\cr\end{matrix}\right)$ be a matrix in $SL(2,\bZ) $ such that $c \xi  +d > 0$. Then the point
$
\left(\begin{matrix}\xi' \cr y'\cr\end{matrix}\right)
$
with coordinates
$$
\xi' = {a\xi + b\over  c\xi +d} \quad {\rm and}\quad y' = {y\over c\xi +d}
$$
belongs to $\cE'(\psi)$. Indeed, substituting 
 $$
 q = a q' +cp' , \quad p = bq' +d p' \leqno{(10)}
 $$
in $(9)$ and dividing by $c\xi +d$, we obtain   the inequalities
 $$
 q' \ge 1 \quad{\rm and} \quad  | q'  \xi'  + p' -y' | \le {\kappa\over c\xi +d}  \psi( q) \le  \kappa' \psi(q') , \leqno{(11)}
 $$
 for some $\kappa' >0$ independent of $q'$. The positivity of $q'$ is proved  as follows. Note that $(9)$ implies
the estimate
$$
p = -q \xi +\cO_{\xi,y}(1).
$$
Then, inverting the linear substitution $(10)$, we find  
 $$
 q' = d q- cp  = q(c\xi + d) + \cO_{\gamma,\xi,  y}(1).
 $$
 Since we have assumed that $c\xi +d >0$, the term $q(c\xi + d)$  is arbitrarily large  when $q$ is large enough. 
The  conditions $(3)$ now show that $\psi(q)\asymp \psi(q')$.
Thus $(11)$ is satisfied for infinitely many primitive points $(p',q')$, since the  linear substitution $(10)$ is unimodular.
  We have shown  that  $\left(\begin{matrix}\xi' \cr y'\cr\end{matrix}\right)$ belongs to $\cE'(\psi)$.

We now  prove  that the intersection $\cE'(\psi) \cap ( \bR \times \bR^+)$ is either a full  or a null subset of the half plane  $\bR \times \bR^+$.
To that purpose,  we consider  the map
$$
\Phi : \bR \times \bR^+ \rightarrow \bR \times \bR^+ ,\quad{ \rm defined \,\,\, by } \quad  \Phi\left(\left(\begin{matrix}x \cr y\cr\end{matrix}\right)\right) 
=\left(\begin{matrix}x/y \cr 1/y\cr\end{matrix}\right).
$$
 Clearly  $\Phi$ is  a continuous involution  of 
$\bR \times \bR^+$. The image 
$$
\Omega : = \Phi \Big(\cE'(\psi ) \cap ( \bR\times \bR^+)\Big)
$$
is formed   by all points of the type
$$
\left(\begin{matrix}u\cr v\cr\end{matrix}\right)
  = \left(\begin{matrix}{\xi\over y} \cr {1\over y}\cr\end{matrix}\right),
$$
where $\left(\begin{matrix}\xi \cr y\cr\end{matrix}\right)$ ranges over $\cE'(\psi ) \cap ( \bR\times \bR^+)$. Now, the above condition $c\xi + d >0 $ is obviously equivalent to
$cu+dv >0$ since $y$ is positive. Then, the point 
$$
 \Phi \left( \begin{matrix} au+bv \cr cu+dv\cr\end{matrix} \right) = \left(  \begin{matrix} {au+bv\over cu+dv} \cr  {1\over cu+dv}\cr\end{matrix} \right) 
 =  \left(\begin{matrix}{a\xi + b\over c\xi +d}\cr {y\over c\xi +d}\cr\end{matrix} \right)
 $$
belongs to $\cE'(\psi ) \cap ( \bR\times \bR^+)$,  by the preceding observation. Applying the involution $\Phi$, we find that
$$
\Phi  \left(    \left(\begin{matrix}{a\xi + b\over c\xi +d}\cr {y\over c\xi +d}\cr \end{matrix} \right)   \right) = 
\left(\begin{matrix} au+bv\cr cu+dv\cr \end{matrix}\right)
= \left(\begin{matrix} a& b\cr c& d\cr \end{matrix}\right)\left(\begin{matrix} u\cr v\cr\end{matrix}\right)
$$
belongs to $\Omega$. In other words, setting   $\Gamma = SL(2,\bZ)$, we have established the inclusion
$$
(\Gamma \Omega) \cap ( \bR\times \bR^+) \subseteq \Omega.
$$
Since the reversed inclusion is obvious, the equality 
$
\Omega = (\Gamma \Omega) \cap ( \bR\times \bR^+) 
$
holds in fact. 
Assuming that $\Omega$ is not a null set, the ergodicity of the linear action of $\Gamma$ on $\bR^2$  \cite{Mo}
shows that $\Gamma \Omega$ is a full set in $\bR^2$. Hence  $\Omega $
 is a full set in the half plane
$\bR\times \bR^+$.  Transforming now $\Omega$  by $\Phi$, we find that 
$$
\Phi(\Omega) = \cE'(\psi) \cap ( \bR \times \bR^+),
$$
is as well  a full set in $\bR\times \bR^+$,  thus proving the claim.

We finally  use another transformation   to carry  the zero-one law  from the positive  half plane $\bR\times \bR^+$ to  
the negative one $\bR\times \bR^-$.
Writing $(9)$ in the equivalent form
$$
q \ge 1 \quad {\rm and}\quad | q(-\xi) + (-p) -( -y) |  \le \kappa \psi( q), 
$$
shows that $\cE'(\psi)$ is invariant under the symmetry $ \left( \begin{matrix} \xi \cr y\cr \end{matrix}\right) \mapsto  \left( \begin{matrix}- \xi \cr - y\cr \end{matrix}\right) $
 which maps 
$\bR\times \bR^+$ onto  $\bR\times \bR^-$. 
Therefore $\cE'(\psi) \cap(\bR\times \bR^-)$
is a null or a full set in $\bR\times \bR^-$ when $\cE'(\psi) \cap(\bR\times \bR^+)$
is accordingly a null or a full set in $\bR\times \bR^+$.
\cqfd

\bigskip
Now, the combination of Lemma 3 and Lemma 4 obviously yields our proposition. 

  \subsection{Concluding the proof of Theorem 2}

   Assume first that $\sum\psi(\ell)$ converges. We have to show that the set
  $$
  \cE(\psi) = \limsup_{q \rightarrow + \infty}\cE_q(\psi)
  $$
  has null Lebesgue measure. Lemma 2  shows that the partial sums
  $$
  \sum_{q=1}^Q \lambda( \cE_q(\psi) \cap [0,1]^2) = 2 \sum_{q=1}^Q{\varphi(q)\psi(q)\over q} \le 2 \sum_{q=1}^Q\psi(q)
  $$
  converge (*).
  \footnote{(*) Here again we assume without loss of generality that $\psi(q) \le 1/2$ for every $q\ge 1$, so that Lemma 2 may be applied.}
  Then, Borel-Cantelli Lemma ensures that the $\limsup$  set  $\cE(\psi)\cap[0,1]^2$ is a null set.
  Thus $\cE(\psi)$ cannot be a full set. Now, the above proposition tells us that $\cE(\psi)$ is a null set.

   We now consider  the case of a divergent series $\sum\psi(\ell)$. 
  Observe  that the estimate
  $$
  {1\over 2} \sum_{q=1}^Q\psi(q) \le \sum_{q=1}^Q {\varphi(q) \psi (q)\over q} \le \sum_{q=1}^Q \psi(q)
  \leqno{(12)}
  $$
  holds true for any large integer $Q$, since the sequence $\psi(\ell)_{\ell\ge 1}$ is non-increasing. The right inequality is obvious,
  while the left one  easily follows from  Abel  summation process. See for instance Chapter 2 of \cite{Har}, where full details are provided.
  By Lemma 2 and $(12)$, the sums
  $$
 \sum_{q=1}^Q \lambda( \cE_q(\psi) \cap [0,1]^2) = 2 \sum_{q=1}^Q{\varphi(q)\psi(q)\over q} \ge   \sum_{q=1}^Q\psi(q)
  $$ 
  are then unbounded. Then, using a   classical  converse to  Borel-Cantelli Lemma,  we have   the  lower bound 
 $$
 \begin{aligned}
   \lambda\Big( \cE(\psi) \cap [0,1]^2\Big) &=  \lambda\Big( \limsup_{q \rightarrow + \infty}(  \cE_q(\psi) \cap [0,1]^2)\Big)
   \\
  & \ge \limsup_{Q \rightarrow + \infty}{ \left( \sum_{q=1}^Q \lambda (  \cE_q(\psi) \cap [0,1]^2) \right)^2
   \over
   \sum_{q=1}^Q\sum_{s=1}^Q  \lambda (  \cE_q(\psi)\cap \cE_s(\psi) \cap [0,1]^2) }.
   \end{aligned}
   \leqno{(13)}
   $$
    See for instance  Lemma 2.3 in \cite{Har}.  Lemma 2 and $(12)$ now show
   that the numerator on  the right  hand side of $(13)$  equals
  $$
  4  \left( \sum_{q=1}^Q {\varphi(q)\psi(q)\over q} \right)^2
   \ge 
 \left( \sum_{q=1}^Q \psi(q) \right)^2,
  $$
  when $Q$ is large, while the denominator is bounded from above by
  $$
  4 \displaystyle \sum_{{q=1,  s=1\atop q\not= s}}^Q \psi(q)\psi(s) + 2 \sum_{q=1}^Q \psi(q) 
  \le 4  \left( \sum_{q=1}^Q \psi(q) \right)^2 +  2 \sum_{q=1}^Q \psi(q).
  $$
  Thus $(13)$ yields the lower bound
  $$
   \lambda\Big( \cE(\psi) \cap [0,1]^2\Big) \ge {1\over 4}.
   $$
  Hence $\cE(\psi)$ is not a null set;    it is thus   a full set according to our proposition.

    \section{ An  approach to our problem}

  In this section, we apply  a  transference principle between homogeneous and inhomogeneous 
  approximation, as displayed  in  Chapter V of \cite{Cas} and in \cite{BuLa}, for   constructing  explicit integer solutions of  the inequality
  $$
   | q \xi  + p -y | \le { 2 \over | q |  }. \leqno{(14)}
   $$
   
   Let  $(p_k/q_k)_{k\ge 0}$ be the sequence of convergents to   the irrational number $\xi$. 
     The theory of continued fractions, see for instance the monograph \cite{Kh},
   tells us that 
  $$
   | q_k \xi -p_k | \le {1 \over   q_{k+1}} \quad {\rm and}\quad  p_kq_{k+1} - p_{k+1}q_k = (-1)^{k+1},  \leqno{(15)}
  $$
  for any $k\ge 0$.
Setting 
$
\nu_k = (-1)^{k+1}q_k y
$, we thus  have the relations
$$
\nu_k q_{k+1} +\nu_{k+1} q_k  =0
\quad {\rm and}\quad
\nu_k(q_{k+1} \xi -p_{k+1} ) + \nu_{k+1} (q_k\xi -p_k) = y . \leqno{(16)}
$$
Now, let $n_k$  be anyone of the two integers $\lfloor\nu_k\rfloor$ and  $\lceil\nu_k\rceil$ (\dag).
   \footnote{ (\dag) As usual  $\lfloor x \rfloor$ and $\lceil x \rceil$ stand  respectively for the integer part and the upper integer part of 
  the real number $x$. Then $\lceil x \rceil = \lfloor x \rfloor + 1$, unless $x$ is an integer in which case $\lfloor x \rfloor = \lceil x \rceil = x$.}
  Then,
    $$
  | \nu_k - n_k |  <    1, \leqno{(17)}
  $$
  and $n_k$ is either equal to $(-1)^{k+1}\lfloor yq_k\rfloor$ or to $(-1)^{k+1}\lceil yq_k\rceil$. 
 Setting
  $$
  p = - n_k p_{k+1} -  n_{k+1} p_k \quad {\rm and}\quad q = n_k q_{k+1} + n_{k+1} q_k,  \leqno{(18)}
$$
  we deduce from $(16)$ the expressions  
  $$
  \begin{aligned}
    q  \xi  + p -y     & =    n_k(q_{k+1} \xi -p_{k+1} ) + n_{k+1} (q_k\xi -p_k)-y  
   \\
   & =  (n_k-\nu_k)(q_{k+1} \xi -p_{k+1} ) +(n_{k+1}- \nu_{k+1}) (q_k\xi -p_k)
    \end{aligned}
  \leqno{(19)}
  $$
   and
   $$
    q     =  (n_k- \nu_{k})q_{k+1} + (n_{k+1} -\nu_{k+1})q_k  . 
   \leqno{(20)}
  $$
  Recall that $q_k\xi -p_k$ and $q_{k+1}\xi - p_{k+1}$ have opposite signs.  Assuming that $n_k-\nu_k$ 
 and $n_{k+1} -\nu_{k+1}$ have the same  sign, we infer from  the formulas $(19), (20)$ and from $(15), (17)$ that
 $$
   | q  \xi  + p -y  | < { 1 \over q_{k+1}} \quad { \rm and} \quad | q | <  2q_{k+1}. 
   \leqno{(21)}
   $$
   Otherwise, we have
  $$
   | q  \xi  + p -y  | < { 2 \over q_{k+1}} \quad { \rm and} \quad | q | <  q_{k+1}. 
   \leqno{(22)}
   $$
The  inequalities $(21)$ and $(22)$  obviously imply $(14)$. 
  
  Since the linear substitution $(18)$ is unimodular, the integers $p$ and $q$ are coprime if and only if 
  $n_k$ and $n_{k+1}$ are  coprime. Recall  that   the two  choices $n_k =\lfloor \nu_k \rfloor$
  and $n_k = \lceil \nu_k \rceil $ are admissible, both for $n_k$ and $n_{k+1}$.
  It thus remains to find  indices $k$ for which  at least one 
  of the  coprimality conditions 
  $$
  \begin{aligned}
& \gcd(\lfloor y q_k \rfloor, \lfloor y q_{k+1} \rfloor) = 1   \quad{\rm or} \quad     \gcd(\lceil y q_k \rceil , \lceil  y q_{k+1} \rceil) =1  
\\
 \quad {\rm or}\quad  &  \gcd(\lfloor y q_k \rfloor, \lceil  y q_{k+1} \rceil) =1 \quad{\rm or} \quad   
       \gcd(\lceil y q_k \rceil , \lfloor  y q_{k+1} \rfloor) =1, 
       \end{aligned}
       \leqno{(23)}
  $$
  is verified.  Note that  $(23)$ obviously  fails for all  $k \ge 0$ when  $y$ is an  integer not equal to $1$ or to $-1$.
  Otherwise, the contingent existence of  infinitely many  indices $k$ satisfying $(23)$ is a non-trivial problem   that we leave hanging.
  
  \subsection
  {Proof of Theorem 3}  
    We quote  the following metrical result due to Harman (Theorem 8.3 in \cite{Har}). Assume that
  the series $(5)$
   diverges. Then for almost all positive real numbers $y$, there exist infinitely many indices $k$ such that the integer
  part $\lfloor yq_k \rfloor $ is a prime number. These indices $k$  fulfill   $(23)$ since, assuming for simplicity that $y$ is irrational,
   either $\lfloor yq_{k+1}\rfloor$
  or $\lceil y q_{k+1} \rceil = \lfloor yq_{k+1}\rfloor +1$ is not divisible by   $\lfloor yq_k \rfloor $ and is thus  relatively prime with $\lfloor yq_k \rfloor $. 
  Hence $(14)$ has infinitely many coprime solutions $(p,q)$ for almost 
  every positive real number $y$.  Writing  now $(14)$ in the equivalent form 
  $$
  | (-q) \xi +(-p) - (-y) | \le {2 \over | q |}
  $$
  shows that, $\xi$ being given,  the set of all real numbers $y$ for which $(14)$ has infinitely many coprime solutions is invariant by the symmetry $y\mapsto -y$.
 The  first assertion is thus established. To complete the proof, note that 
  $$ 
  \lim_{k \rightarrow +\infty} { \log q_k \over k} = { \pi^2\over 12 \log 2}
  $$
  for almost every $\xi$ by Khintchine-Levy  Theorem (see equation (4.18) in \cite{Bi}). Thus the series $(5)$ diverges for almost every $\xi$.

 % \medskip
   % Our  problem can be rephrased as finding a primitive point inside thin  strips  on  the plane. 
    %In another direction  let us finally  mention the work  \cite{No},  where are obtained 
  % asymptotic estimates for  the number of bounded primitive points located  in a  strip centered at the origin. 
   % Theorem 3.1 of \cite{No} provides such an estimate  formulated in terms of Farey sequences.
   % Counting formulas could hopefully be employed for improving Theorem 1. 
      
    \section{ Generic density exponents}

We prove in this section  Theorem 4, as a consequence of Borel-Cantelli Lemma combined with the following counting result.

\begin{lem}
Let $\bx$ be a point in $ \bR^2$ whose  orbit $\Gamma \bx$ is dense in $\bR^2$.  For every symmetric 
compact set $\Omega$ in $\bR^2\setminus\{\bf 0\}$ there exists $c > 0$  such that 
$$
 \Card \{ \gamma \in \Gamma ; \gamma \bx \in \Omega , | \gamma | \le T \} \le c T
$$
for any real number $T \ge 1$.
\end{lem}

\pro
Ledrappier \cite{Le} has shown that  the  limit  formula
$$
\lim_{T \rightarrow + \infty}{1\over T} {\displaystyle \sum_{\gamma\in \Gamma, | \gamma | \le T}f(\gamma \bx)}
= {4 \over | \bx | \vol(\Gamma \setminus SL(2,\bR))}\int  { f(\by)  \over | \by |} d\by, 
$$
holds  for any  even continuous function $f : \bR^2\rightarrow \bR$ having compact support on $\bR^2\setminus\{\bf 0\}$, 
with   a suitable normalisation of Haar  measure on $SL(2,\bR)$.
Approximating uniformly from above and from below the characteristic function of $\Omega$ by even continuous functions, we deduce that
$$ 
\lim_{T \rightarrow + \infty}{\Card \{ \gamma \in \Gamma ; \gamma \bx \in \Omega , | \gamma | \le T \}\over T}
= {4 \over | \bx | \vol(\Gamma \setminus SL(2,\bR))}\int _\Omega { d\by \over | \by |}. 
$$
 Lemma 5 immediately follows. \cqfd

For any point $\by\in \bR^2$ and any positive real number $r$, we denote by
$$
B(\by,r) = \{ \bz \in \bR^2 ; | \bz - \by | \le r\}
$$
the closed disc centered at $\by$ with radius $r$. 

\begin{lem}
Let $\bx$ be a point in $ \bR^2$ whose  orbit $\Gamma \bx$ is dense, 
 $\Omega$  a symmetric  compact set  in $\bR^2\setminus\{\bf 0\}$
and   $\mu$  a real number $> 1/2$. 
For  every  integer $n \ge 1$, put
$$
\cB_n = \bigcup_{{\gamma \in \Gamma \atop | \gamma | = n, \gamma \bx \in \Omega}} B(\gamma \bx , n^{-\mu}).
$$
Then the set
$$
\cB  : = \limsup_{n \rightarrow +\infty} \cB_n = \bigcap_{N\ge 1} \bigcup_{n \ge N}\cB_n = 
 \bigcap_{N\ge 1} \bigcup_{{\gamma \in \Gamma \atop | \gamma | \ge N, \gamma \bx \in \Omega}}B( \gamma \bx, | \gamma |^{-\mu})
 $$
 has null Lebesgue measure.
\end{lem}

\pro 
 We apply Borel-Cantelli Lemma and we prove 
that the series  $\sum_{n\ge 1}\lambda(\cB_n)$ converges if $\mu >1/2$.

 For every  positive integer $n$, set
$$
M_n =  \Card \{ \gamma \in \Gamma ; \gamma \bx \in \Omega, | \gamma | = n \}.
$$
 Lemma 5 gives us   the upper bound 
$$
M_1 + \cdots + M_n =  \Card \{ \gamma \in \Gamma ; \gamma \bx \in \Omega, | \gamma | \le n \} \le c n  ,  \leqno{(24)}
$$
 for some $c>0$ independent of $n \ge 1$.
Since a ball of radius $r$ has Lebesgue measure $4r^2$, we trivially bound from above
$$
\lambda( \cB_n) \le  \sum_{{\gamma \in \Gamma \atop | \gamma | = n, \gamma \bx \in \Omega}}4n^{-2\mu} = 4  M_n n^{-2\mu}.
$$
Summing by parts, we deduce from $(24)$ that 
$$
\begin{aligned}
\sum_{n= 1}^N {M_n \over n^{2\mu}}  & = \sum_{n= 1}^{N-1} (M_1 + \cdots + M_n)\left( {1\over n^{2\mu}}- {1\over (n+1)^{2\mu}}\right) 
+ {M_1+ \cdots + M_N\over N^{2\mu}}
\\
& \le c   \sum_{n=1}^{N-1} n \left( {1\over n^{2\mu}} - {1\over (n+1)^{2\mu}}\right)  +  {c  N \over N^{2\mu}} = c  \sum_{n=1}^{N}  {1\over n^{2\mu}}.
\end{aligned}
$$
The partial sums 
$$
\sum_{n=1}^N \lambda(\cB_n) \le 4 \sum_{n= 1}^N {M_n \over n^{2\mu}} \le 4 c \sum_{n=1}^{N}  {1\over n^{2\mu}}
$$
thus converge  if $\mu > 1/2$. \cqfd

\subsection
{Proof of Theorem 4}
We argue by contradiction and suppose on the contrary that $\mu_\Gamma(\bx) > 1/2$.
Fix a real number $\mu$ with $1/2 < \mu < \mu_\Gamma(\bx)$. Then for almost all points  $\by\in \bR^2$, we have 
$\mu(\bx,\by) > \mu$. This means that there exist infinitely many $\gamma \in \Gamma$ satisfying (6),  or equivalently that
$\by $ belongs to infinitely many balls of the form $B(\gamma \bx, | \gamma |^{-\mu})$. 
We now restrict our attention to points $\by$
with $\mu(\bx,\by) > \mu$ lying in an annulus 
$$
\Omega' = \{ \bz \in \bR^2 ; a' \le | \bz | \le b' \},
$$
where  $b'>a'>0$ are arbitrarily fixed. Since   $\by$ belongs to 
the intersection 
$
\Omega' \cap B(\gamma \bx, | \gamma |^{-\mu}), 
$
 we deduce  from the triangle inequality the estimate
$$
a'-|\gamma|^{-\mu} \le | \gamma \bx | \le b' + |\gamma|^{-\mu}. 
$$
Fixing $a< a'$ and $b> b'$, the center $\gamma \bx$ then lies in the larger annulus 
$$
\Omega = \{ \bz \in \bR^2 ; a \le | \bz | \le b \},
$$
provided that $|\gamma| $ is large enough. It follows that $\by$ falls inside  the union of balls 
$$
 \bigcup_{{\gamma \in \Gamma \atop | \gamma | \ge N, \gamma \bx \in \Omega}}B( \gamma \bx, | \gamma |^{-\mu})
$$
considered in Lemma 6 for every integer    $N$ large enough, and thus $\by$ belongs to $\cB$. However, Lemma 6 asserts that $\cB$ is a null set
which is a contradiction.

    %As a related result, let us mention the paper \cite{Wa}, where Watson gives the density of yhe  integers $k$ satisfying 
 % $$
  %1 \le k \le n \quad {\rm and}\quad \gcd( k, \lfloor y k \rfloor ) =1,
  %$$
  %as $n$ tends to infinity.

\vskip1cm
 Institut de Math\'ematiques de Luminy,  Case 907, 163 avenue de Luminy, 13288,  Marseille C\'edex 9.

  {\tt michel-julien.laurent@univmed.fr}, 
  
  {\tt  arnaldo.nogueira@univmed.fr}

\end{document}